\documentclass[a4paper,12pt]{article}
\usepackage{amsmath}
\usepackage{amsthm}
\usepackage{amssymb}
\usepackage{enumerate}
\usepackage{url}
\usepackage[utf8]{inputenc}

\title{Pitman's and L\'evy's theorems for Brownian bridges}

\author{Yuu Hariya}

\date{\empty}
\numberwithin{equation}{section}

\theoremstyle{plain}

\newtheorem{theorem}{Theorem}[section]
\newtheorem{lemma}[theorem]{Lemma}

\newtheorem{corollary}[theorem]{Corollary}

\theoremstyle{definition}

\theoremstyle{remark}
\newtheorem{remark}[theorem]{Remark}

\begin{document}

\newcommand\ND{\newcommand}
\newcommand\RD{\renewcommand}

\ND\R{\mathbb{R}}

\ND\ind{\boldsymbol{1}}

\ND\ga{\gamma}
\ND\la{\lambda }
\ND\ve{\varepsilon}

\ND\lref[1]{Lemma~\ref{#1}}
\ND\tref[1]{Theorem~\ref{#1}}
\ND\pref[1]{Proposition~\ref{#1}}
\ND\sref[1]{Section~\ref{#1}}
\ND\ssref[1]{Subsection~\ref{#1}}
\ND\aref[1]{Appendix~\ref{#1}}
\ND\rref[1]{Remark~\ref{#1}} 
\ND\cref[1]{Corollary~\ref{#1}}
\ND\csref[1]{Corollaries~\ref{#1}}
\ND\eref[1]{Example~\ref{#1}}
\ND\fref[1]{Fig.\ {#1} }
\ND\lsref[1]{Lemmas~\ref{#1}}
\ND\tsref[1]{Theorems~\ref{#1}}
\ND\dref[1]{Definition~\ref{#1}}
\ND\psref[1]{Propositions~\ref{#1}}
\ND\rsref[1]{Remarks~\ref{#1}}
\ND\sssref[1]{Subsections~\ref{#1}}

\ND\pr{\mathbb{P}}
\ND\ex{\mathbb{E}}

\ND\eqd{\stackrel{(d)}{=}}

\ND\cm{\mathcal{M}}

\ND\cp{\mathcal{P}}

\ND\sg{\sigma}
\ND\ctd[2]{C([0,#1];\mathbb{R}^{#2})}

\ND\pbes[2]{\mathbb{P}^{(3)}_{#1,#2}}
\ND\ebes[2]{\mathbb{E}^{(3)}_{#1,#2}}

\ND\bb{\beta}

\ND\bd{\mathbf{B}}

\ND\cl{\mathcal{L}}

\ND{\rmid}[1]{\mathrel{}\middle#1\mathrel{}}

\def\thefootnote{{}}

\maketitle 
\begin{abstract}
A synthetic study of Pitman's and L\'evy's theorems for one-dimensional Brownian bridges with arbitrary endpoints is provided. 
\footnote{Mathematical Institute, Tohoku University, Aoba-ku, Sendai 980-8578, Japan}
\footnote{E-mail: hariya@tohoku.ac.jp}
\end{abstract}

\section{Introduction and main results}\label{;intro}

Let $B=\{ B_{s}\} _{s\ge 0}$ be a standard one-dimensional Brownian motion 
and $R=\{ R_{s}\} _{s\ge 0}$ a three-dimensional Bessel process starting from 
the origin. On the space $C([0,\infty );\R )$ of real-valued continuous functions 
$\phi _{s},s\ge 0$, over the half-line $[0,\infty )$, we define the transformation $\cp $ by 
\begin{align*}
 \cp (\phi )(s)&:=2\max _{0\le u\le s}\phi _{u}-\phi _{s},\quad s\ge 0. 
\end{align*}
The celebrated Pitman's theorem \cite{jwp} asserts that the process 
$\{ \cp (B)(s)\} _{s\ge 0}$ has the same law as $R$, or more precisely 
(see Theorem~1.3 in \cite{jwp} and the description of the theorem that 
follows its statement), 
\begin{align}\label{;pit}
 \bigl\{ 
 \bigl( \cp (B)(s),B_{s}\bigr) 
 \bigr\} _{s\ge 0}\eqd 
 \Bigl\{ 
 \Bigl( 
 R_{s},-R_{s}+2\inf _{u\ge s}R_{u}
 \Bigr) 
 \Bigr\} _{s\ge 0}.  
\end{align}
Here the equality stands for the identity in law. See also 
\cite[Chapter~VI, Theorem~\thetag{3.5}]{ry} for the above formulation 
of the theorem.   

Let $L=\{ L_{s}\} _{s\ge 0}$ be the local time process at $0$ of $B$. 
L\'evy's theorem tells us that 
\begin{align}\label{;lev}
\left\{ 
\bigl( |B_{s}|+L_{s},|B_{s}|\bigr) 
\right\} _{s\ge 0}
\eqd 
\Bigl\{ 
\Bigl( 
\cp (B)(s),\max _{0\le u\le s}B_{u}-B_{s}
\Bigr) 
\Bigr\} _{s\ge 0}
\end{align}
(see, e.g., \cite[Chapter~VI, Theorem~\thetag{2.3}]{ry}). Both of the 
two theorems have attracted many researchers and a number of 
studies have been conducted such as extensions to Brownian motion 
with drift, to geometric Brownian motions, to the solutions of a 
class of stochastic differential equations driven by Brownian motion, 
to L\'evy processes, and so forth; see 
\cite{ber, mySII, rv, rvy, st} and references therein. 

Of concern in this paper is a natural fundamental question of what if 
the underlying process is a Brownian bridge. In what follows, unless 
otherwise specified, we fix $t>0$, and for every $x\in \R $, we let 
$\bb ^{x}=\{ \bb ^{x}_{s}\} _{0\le s \le t}$ be a one-dimensional Brownian 
bridge from $0$ to $x$ over the time interval $[0,t]$. We recall that 
$\bb ^{x}$ is related with $\bb ^{0}$ via 
\begin{align}\label{;bb0}
 \bb ^{x}\eqd \left\{ \bb ^{0}_{s}+\frac{x}{t}s\right\} _{0\le s \le t},
\end{align}
and that the process $\bb ^{0}$ is a continuous centered 
Gaussian process with covariance function 
\begin{align*}
 \frac{s_{1}(t-s_{2})}{t},\quad 0\le s_{1}\le s_{2}\le t.
\end{align*}
See \cite[Section~1.1]{my} for an account of one-dimensional 
Brownian bridges. As for the question raised above, while the case 
$x=0$ is well understood especially in the context of path decomposition 
(see \cite{ber0, bp, by} and references therein), it seems that less is known 
in the case $x\neq 0$. The aim of this paper is to provide a synthetic 
treatment of the question and give a complete answer; we will do this 
by employing \tsref{;tmaini} and \ref{;tmainii} below, whose statements 
are of interest in their own right.

Let $\bd =\bigl\{ \bd _{s}=(B^{1}_{s},B^{2}_{s},B^{3}_{s})\bigr\} _{s\ge 0}$ 
be a standard three-dimensional Brownian motion, namely 
$B^{i}=\{ B^{i}_{s}\} _{s\ge 0},i=1,2,3$, are independent, standard 
one-dimensional Brownian motions. For each point $(x_{1},x_{2},x_{3})$ in $\R ^{3}$, 
we write $|(x_{1},x_{2},x_{3})|$ for its Euclidean norm 
$\sqrt{x_{1}^{2}+x_{2}^{2}+x_{3}^{2}}$. One of the main results of the 
paper is stated as 

\begin{theorem}\label{;tmaini}
The two-dimensional process 
\begin{align*}
 \bigl( \cp (B)(s),B_{s}\bigr) ,\quad 0\le s\le t,
\end{align*}
is identical in law with the process given by 
\begin{align}\label{;qtmaini}
 \Bigl( 
 |\bd _{s}|,-|\bd _{s}|+\min \Bigl\{ 
 2\min _{s\le u\le t}|\bd _{u}|,|\bd _{t}|+B^{1}_{t}
 \Bigr\} 
 \Bigr) ,\quad 0\le s\le t.
\end{align}
\end{theorem}

Given $x\in \R $, supposing that the Brownian bridge $\bb ^{x}$ is 
independent of the two-dimensional Brownian motion $(B^{2},B^{3})$, 
we let $M^{x}=\{ M^{x}_{s}\} _{0\le s\le t}$ be a process that is identical 
in law with the process given by 
\begin{align*}
 \sqrt{(\bb ^{x}_{s})^{2}+(B^{2}_{s})^{2}+(B^{3}_{s})^{2}},\quad 0\le s\le t.
\end{align*}
For every $y\in \R $, we define the transformation $\cl _{y}$ on 
the space $\ctd{t}{}$ of real-valued continuous functions over $[0,t]$ by 
\begin{align*}
 \cl _{y}(\phi )(s):=-\phi _{s}+\min \Bigl\{ 
 2\min _{s\le u\le t}\phi _{u},\phi _{t}+y
 \Bigr\} 
\end{align*} 
for $\phi \in \ctd{t}{}$ and $0\le s\le t$, so that the second component 
in \eqref{;qtmaini} is expressed as $\cl _{B^{1}_{t}}(|\bd |)(s)$ for all $0\le s\le t$. 
Notice that 
\begin{align*}
 \cl _{B^{1}_{t}}(|\bd |)(t)&=-|\bd _{t}|+\min \bigl\{ 2|\bd _{t}|,|\bd _{t}|+B^{1}_{t}\bigr\} \\
 &=-|\bd _{t}|+|\bd _{t}|+B^{1}_{t}\\
 &=B^{1}_{t},
\end{align*}
and hence \tref{;tmaini} entails the joint identity in law 
\begin{align}\label{;jil}
 \Bigl( 
 \bigl\{ 
 \bigl( 
 \cp (B)(s),B_{s}
 \bigr) 
 \bigr\} _{0\le s\le t},B_{t}
 \Bigr) 
 \eqd 
 \Bigl( 
 \bigl\{ 
 \bigl( 
 |\bd _{s}|,\cl _{B^{1}_{t}}(|\bd |)(s)
 \bigr) 
 \bigr\} _{0\le s\le t},B^{1}_{t}
 \Bigr) .
\end{align}
Conditionally on the second component on the respective sides with 
$B_{t}=x$ and $B^{1}_{t}=x$, we immediately obtain the following corollary, 
Pitman's theorem for Brownian bridges. In the statement, we restrict 
the Pitman transformation $\cp $ to $\ctd{t}{}$ keeping the same notation.

\begin{corollary}\label{;cmaini}
For every $x\in \R $, it holds that 
\begin{align}\label{;qcmaini}
\left\{ 
\bigl( 
\cp (\bb ^{x})(s),\bb ^{x}_{s}
\bigr) 
\right\} _{0\le s\le t}
\eqd 
\left\{ 
\bigl( 
M^{x}_{s},\cl _{x}(M^{x})(s)
\bigr) 
\right\} _{0\le s\le t}.
\end{align}
\end{corollary}

The above corollary suggests that 
\begin{align}\label{;scl}
 \cl _{\phi _{t}}(\cp (\phi ))(s)=\phi _{s}
\end{align}
for all $\phi \in \ctd{t}{}$ and $0\le s\le t$, which is 
indeed the case since, in the notation of \cite{har25+} 
by the author, the left-hand side is nothing but 
$\cm _{\phi _{t}}(\phi )(s)$ in view of Proposition~2.1 therein. 
See also \rref{;rscl} below.

We give some remarks on the special case $x=0$.
In what follows, we shall refer to a function 
on $\ctd{t}{}$ valued in $\overline{\R }$, the extended real line,   
as a measurable function on this space if it is Borel-measurable 
with respect to the topology of uniform convergence.

\begin{remark}
\thetag{1} When $t=1$, the process $M^{0}$ is a Brownian meander; 
see \cite[p.~50, Corollary~3.9.1]{my}. (The case of general $t$ is referred 
to as a Brownian meander of duration $t$ \cite{imh}.) We may compare 
the above \cref{;cmaini} with Theorems~2.2 and 4.3 in \cite{bp}, the latter 
of which corresponds to the identity in law between the first components 
in \eqref{;qcmaini} and the former to that between the second 
components in \eqref{;qcmaini}. Our corollary unifies those two 
theorems, with the difference that, by adopting the notation  
\begin{align*}
 U=\sup \left\{ 
 0\le s\le 1;\,M^{0}_{s}=\frac{1}{2}M^{0}_{1}
 \right\} 
\end{align*}
from \cite{bp} (in fact, $U$ is uniformly distributed on $[0,1]$), 
Theorem~2.2 in \cite{bp} asserts that the process 
\begin{align}\label{;utr}
 M^{0}_{U-s}-\frac{1}{2}M^{0}_{1},\quad 0\le s\le U,
\end{align}
is identical in law with $\bb ^{0}$ running until it reaches 
the global minimum, while our \cref{;cmaini} shows that the process 
\begin{align*}
 M^{0}_{s}-2\min _{s\le u\le 1}M^{0}_{u},\quad 0\le s\le U,
\end{align*}
is.\ This difference may be explained by Theorem~1.1 in \cite{har25+}, 
which will be pursued elsewhere. Our corollary provides a more direct 
connection between $\bb ^{0}$ and $M^{0}$ in that it does not 
involve a time-reversal operation such as \eqref{;utr}, which 
presumably has enabled us to treat the laws of $\cp (\bb ^{0})$ 
and $\bb ^{0}$ simultaneously, namely their joint identity in law.

\thetag{2} Imhof's relation \cite[Section~4]{imh} shows the following 
absolute continuity relationship between the laws of $M^{0}$ and 
the three-dimensional Bessel process $R$ restricted to the time 
interval $[0,t]$: for every nonnegative measurable function $F$ on $\ctd{t}{}$,
\begin{align}\label{;imrel}
 \ex \!\left[ 
 F(M^{0})
 \right] 
 =\ex \biggl[ 
 F(R)\sqrt{\frac{\pi t}{2}}\frac{1}{R_{t}}
 \biggr] ;
\end{align}
we also refer to \cite[Th\'eor\`eme~3]{by}. 
\end{remark}

Combining \cref{;cmaini} with Theorem~1.4 in \cite{har25+}, we 
obtain yet another corollary.

\begin{corollary}\label{;ccmaini}
For every $x\in \R $, it holds that 
\begin{align*}
 \ex \!\left[ 
 F(M^{x})
 \right] =
 \ex \!\left[ 
 F(M^{0})\rmid| M^{0}_{t}\ge |x|
 \right] 
\end{align*}
for any nonnegative measurable function $F$ on $\ctd{t}{}$.
\end{corollary}

It follows from the above two corollaries that $\cp (\bb ^{x})$ 
has the same law as $M^{0}$ conditioned on the event that 
$M^{0}_{t}\ge |x|$. This conditional identity in law is observed 
in \cite{bp}, in the remark that follows the description of Theorem~4.3, 
where the above-mentioned conditional law of $M^{0}$ is not 
identified with $M^{x}$.

We turn to the identity~\eqref{;lev} of L\'evy. Similar reasoning to 
the derivation of \tref{;tmaini} enables us to obtain the 
following identity in law.

\begin{theorem}\label{;tmainii}
The two-dimensional process 
\begin{align}\label{;q1tmainii}
 (|B_{s}|+L_{s},|B_{s}|),\quad 0\le s\le t,
\end{align}
is identical in law with the process given by 
\begin{align}\label{;q2tmainii}
 \Bigl( 
 |\bd _{s}|,|\bd _{s}|-\min \Bigl\{ 
 \min _{s\le u\le t}|\bd _{u}|,|\bd _{t}|-|B^{1}_{t}|
 \Bigr\} 
 \Bigr) ,\quad 0\le s\le t.
\end{align}
\end{theorem}

Notice that the evaluation at $s=t$ of the second component in \eqref{;q2tmainii} yields $|B^{1}_{t}|$.
Pick $x\in \R $ arbitrarily. We discuss the conditional law of \eqref{;q1tmainii} given $|B_{t}|=|x|$ and that of 
\eqref{;q2tmainii} given $|B^{1}_{t}|=|x|$. Let $\la (\bb ^{x})=\{ \la _{s}(\bb ^{x})\} _{0\le s\le t}$ be the 
local time process at $0$ of $\bb ^{x}$ that is characterized as 
\begin{align}\label{;char}
 \pr \!\left( 
 \la _{s}(\bb ^{x})=\lim _{\epsilon \downarrow 0}\frac{1}{2\epsilon }\int _{0}^{s}\ind _{(-\epsilon ,\epsilon )}(\bb ^{x}_{u})\,du,\ 
 0\le \forall s\le t
 \right) =1
\end{align}
in view of Theorem~9.4 and Proposition~9.9 in \cite{lg}, as well as the remark just after the 
proposition, due to the fact that 
\begin{align*}
 \int _{0}^{t}\ind _{\{ 0\} }(\bb ^{x}_{s})\left| 
 \frac{x-\bb ^{x}_{s}}{t-s}
 \right| ds=0,\quad \text{a.s.},
\end{align*} 
thanks to Fubini's theorem: 
\begin{align*}
 \ex \!\left[ 
 \int _{0}^{t}\ind _{\{ 0\} }(\bb ^{x}_{s})\left| 
 \frac{x-\bb ^{x}_{s}}{t-s}
 \right| ds
 \right] 
 =\int _{0}^{t}\ex \!\left[ 
 \left| 
 \frac{x-\bb ^{x}_{s}}{t-s}
 \right| ;\bb ^{x}_{s}=0
 \right] ds=0.
\end{align*}
Notice that the process 
\begin{align*}
 \int _{0}^{s}\frac{x-\bb ^{x}_{u}}{t-u}\,du,\quad 0\le s\le t,
\end{align*}
is the finite variation part in the semimartingale decomposition of 
$\bb ^{x}$ in its own filtration.
By relation~\eqref{;bb0} and the symmetry $-\bb ^{0}\eqd \bb ^{0}$, 
we have $-\bb ^{-x}\eqd \bb ^{x}$, from which it is readily seen that 
the law of $\{ B_{s}\} _{0\le s\le t}$ conditioned on $|B_{t}|=|x|$ 
coincides with the law of $\{ \ve \bb ^{x}_{s}\} _{0\le s\le t}$, where 
$\ve $ is a Rademacher (or symmetric Bernoulli) random variable that is 
independent of $\bb ^{x}$. The same observation applies to 
$\{ B^{1}_{s}\} _{0\le s\le t}$. In addition, the above characterization~\eqref{;char} 
entails that the local time process at $0$ of the Brownian bridge $-\bb ^{x}$ 
coincides with that of $\bb ^{x}$:
\begin{align}\label{;chard}
 \la _{s}(-\bb ^{x})=\la _{s}(\bb ^{x}),\quad 0\le s\le t,\ \text{a.s.}
\end{align}
Summing up these observations, we obtain a corollary to 
\tref{;tmainii} as follows:

\begin{corollary}\label{;cmainii}
For every $x\in \R $, the two-dimensional process 
\begin{align}\label{;q1cmainii}
 \left( 
 |\bb ^{x}_{s}|+\la _{s}(\bb ^{x}),|\bb ^{x}_{s}|
 \right) ,\quad 0\le s\le t,
\end{align}
is identical in law with the process given by 
\begin{align}\label{;q2cmainii}
 \Bigl( 
 M^{x}_{s},M^{x}_{s}-\min \Bigl\{ 
 \min _{s\le u\le t}M^{x}_{u},M^{x}_{t}-|x|
 \Bigr\} 
 \Bigr) ,\quad 0\le s\le t.
\end{align}
\end{corollary}

In the case $t=1$ and $x=0$, the above corollary agrees with \cite[Th\'eor\`eme~8]{by}. 
See also \cite[Theorem~3.1]{bp}.

By the fact that 
\begin{align}\label{;eqds}
 M^{x}\eqd M^{-|x|}\eqd \cp (\bb ^{-|x|})
\end{align}
in view of \cref{;ccmaini} and thanks to \cref{;cmaini}, 
we have from \cref{;cmainii} the counterpart of L\'evy's theorem 
in the case of Brownian bridges. For a real number $a$, we write 
$a_{+}=\max \{ a,0\} $.

\begin{corollary}\label{;ccmainii}
For every $x\in \R $, the two two-dimensional processes 
\begin{align*}
 \left( 
 |\bb ^{x}_{s}|+\la _{s}(\bb ^{x}),|\bb ^{x}_{s}|
 \right) ,\quad 0\le s\le t,
\end{align*}
and
\begin{align*}
 \biggl( 
 \cp (\bb ^{-|x|})(s),
 \min \biggl\{ 
 \max _{0\le u\le s}\bb ^{-|x|}_{u},
 \Bigl( 
 \max _{s\le u\le t}\bb ^{-|x|}_{u}
 \Bigr) _{+}
 \biggr\} -\bb ^{-|x|}_{s}
 \biggr) ,\quad 0\le s\le t,
\end{align*}
have the same law.
\end{corollary}

We consider the two mappings $\ga ,\sg \colon \ctd{t}{}\to [0,t]$ defined respectively by 
\begin{align*}
 \ga (\phi ):=\sup \{ 0\le s\le t;\,\phi _{s}=\phi _{0}\} , && 
 \sg (\phi ):=\inf \Bigl\{ 
 0\le s\le t;\,\phi _{s}=\max _{0\le u\le t}\phi _{u}
 \Bigr\} ,
\end{align*}
for $\phi \in \ctd{t}{}$. Observe that, a.s., $\sg (\bb ^{-|x|})\le \ga (\bb ^{-|x|})$ and 
\begin{align*}
  \Bigl( 
 \max _{s\le u\le t}\bb ^{-|x|}_{u}
 \Bigr) _{+}=
 \begin{cases}
 \displaystyle \max _{0\le u\le t}\bb ^{-|x|}_{u}, & 0\le s\le \sg (\bb ^{-|x|}),\\
 \displaystyle \max _{s\le u\le t}\bb ^{-|x|}_{u}, & \sg (\bb ^{-|x|})\le s\le \ga (\bb ^{-|x|}),\\
 0, & \ga (\bb ^{-|x|})\le s\le t,
 \end{cases} 
\end{align*}
since $\max _{0\le u\le t}\bb ^{-|x|}_{u}\ge \bb ^{-|x|}_{0}=0$ and 
$\max _{\ga (\bb ^{-|x|})\le u\le t}\bb ^{-|x|}_{u}=0$, a.s. 
Consequently, the two-dimensional process 
\begin{align*}
 \left( 
 \la _{s}(\bb ^{x}),|\bb ^{x}_{s}|
 \right) ,\quad 0\le s\le t,
\end{align*}
is identical in law with the process described by 
\begin{align*}
 \begin{cases}
 \displaystyle \Bigl( \max _{0\le u\le s}\bb ^{-|x|}_{u}, \max _{0\le u\le s}\bb ^{-|x|}_{u}-\bb ^{-|x|}_{s}\Bigr) , 
 & 0\le s\le \sg (\bb ^{-|x|}),\\
 \displaystyle \Bigl( 2\max _{0\le u\le t}\bb ^{-|x|}_{u}-\max _{s\le u\le t}\bb ^{-|x|}_{u},\max _{s\le u\le t}\bb ^{-|x|}_{u}-\bb ^{-|x|}_{s}\Bigr) , 
 & \sg (\bb ^{-|x|})\le s\le \ga (\bb ^{-|x|}),\\
 \displaystyle \Bigl( 2\max _{0\le u\le t}\bb ^{-|x|}_{u},-\bb ^{-|x|}_{s}\Bigr) ,  
 & \ga (\bb ^{-|x|})\le s\le t.
 \end{cases} 
\end{align*}
The case $t=1$ and $x=0$ recovers \cite[Theorem~4.1]{bp}. If we define 
\begin{align*}
 \tau ^{x}:=\inf \left\{ 
 0\le s\le t;\,\la _{s}(\bb ^{x})=\frac{1}{2}\la _{t}(\bb ^{x})
 \right\} ,
\end{align*}
then we have an identity in law of interest: 
\begin{align*}
 \left( 
 \tau ^{x},\ga (\bb ^{x})
 \right) \eqd 
 \bigl( 
 \sg (\bb ^{-|x|}),\ga (\bb ^{-|x|})
 \bigr) .
\end{align*}
Other interesting distributional identities may also be inferred from \cref{;ccmainii}.

The rest of the paper is organized as follows: in \sref{;sptmains}, we prove \tsref{;tmaini} and \ref{;tmainii};  
in \sref{;spcs}, we prove the four \csref{;cmaini}, \ref{;ccmaini}, \ref{;cmainii} and \ref{;ccmainii} 
to those two theorems.

\section{Proofs of \tsref{;tmaini} and \ref{;tmainii}}\label{;sptmains}

In this section, we prove \tsref{;tmaini} and \ref{;tmainii}. To begin with,  
let $U$ be a random variable subject to the uniform distribution on $[0,1]$ and 
$V$ a random variable identical in law with $2U-1$. Then, $V$ is distributed 
uniformly on $[-1,1]$ and satisfies 
\begin{align}\label{;absv}
 |V|\eqd U.
\end{align}
\begin{lemma}\label{;lu}
It holds that 
\begin{align}
 \Bigl( 
 \{ R_{s}\} _{0\le s\le t},\inf _{u\ge t}R_{u}
 \Bigr) 
 &\eqd \left( 
 \{ R_{s}\} _{0\le s\le t},UR_{t}
 \right) \label{;q1lu}\\
 &\eqd 
 \left( 
 \{ R_{s}\} _{0\le s\le t},(1-U)R_{t}
 \right) ,\label{;q2lu}
\end{align}
where $U$ is assumed to be independent of the three-dimensional 
Bessel process $R$.
\end{lemma}

\begin{proof}
Recall that the global infimum of a three-dimensional Bessel process 
starting from $a>0$ is uniformly distributed on $[0,a]$ (see, e.g., 
\cite[Chapter~VI, Corollary~\thetag{3.4}]{ry}), and hence is identical 
in law with $aU$. Then the Markov property of $R$ entails \eqref{;q1lu}. 
Subsequently, the independence of $R$ and $U$ and the fact that 
$U\eqd 1-U$ prove \eqref{;q2lu}.
\end{proof}

\begin{lemma}\label{;lv}
It holds that 
\begin{align}\label{;qlv}
 \left( 
 \{ R_{s}\} _{0\le s\le t},VR_{t}
 \right) \eqd 
 \left( 
 \{ |\bd _{s}|\} _{0\le s\le t},B^{1}_{t}
 \right) ,
\end{align}
where $V$ is assumed to be independent of $R$. 
\end{lemma}

\begin{proof}
Let $F$ be a nonnegative measurable function on $\ctd{t}{}$ and 
$f$ a nonnegative measurable function on $\R $. What to prove is that 
\begin{align}\label{;q1plv}
 \ex \!\left[ 
 F(|\bd |)f\!\left( 
 \frac{B^{1}_{t}}{|\bd _{t}|}
 \right) 
 \right] =\ex [F(R)]\ex [f(V)].
\end{align}
For every $a\ge 0$, we let $\pbes{0}{a}$ be the law of a three-dimensional 
Bessel bridge from $0$ to $a$ over $[0,t]$ and write $\ebes{0}{a}$ for the 
expectation with respect to it. In view of the rotational invariance of 
multidimensional Brownian motion, we see that, conditionally on 
$\bd _{t}=(x_{1},x_{2},x_{3})\in \R ^{3}$, the law of $\{ |\bd _{s}|\} _{0\le s\le t}$ is 
$\pbes{0}{|(x_{1},x_{2},x_{3})|}$, which may also be deduced from the two general 
facts \cite[Chapter~XI, Exercises~\thetag{3.6} and \thetag{3.7}]{ry} on squared Bessel bridges 
with positive dimensions. Thanks to the above conditional equivalence, 
the left-hand side of the claimed relation~\eqref{;q1plv} is written as 
\begin{align*}
 \int _{\R ^{3}}\ebes{0}{|(x_{1},x_{2},x_{3})|}[F(X)]
 f\!\left( 
 \frac{x_{1}}{|(x_{1},x_{2},x_{3})|}
 \right) 
 \frac{1}{(2\pi t)^{3/2}}
 \exp \left( 
 -\frac{x_{1}^{2}+x_{2}^{2}+x_{3}^{2}}{2t}
 \right) dx_{1}dx_{2}dx_{3},
\end{align*}
where $X=\{ X_{s}\} _{0\le s\le t}$ denotes the coordinate process in $\ctd{t}{}$: 
\begin{align*}
 X_{s}(\phi ):=\phi _{s},\quad 0\le s\le t,\ \phi \in \ctd{t}{}.
\end{align*}
By changing the variables with 
\begin{align*}
 x_{1}=r\cos \theta , && x_{2}=r\sin \theta \cos \varphi ,&& x_{3}=r\sin \theta \sin \varphi ,
\end{align*}
for $r>0$, $0\le \theta \le \pi $ and $0\le \varphi <2\pi $, the above integral is rewritten as 
\begin{align*}
 \int _{0}^{\infty }\ebes{0}{r}[F(X)]\frac{1}{\sqrt{2\pi t^{3}}}r^{2}\exp \left( -\frac{r^{2}}{2t}\right) dr
 \int _{0}^{\pi }f(\cos \theta )\sin \theta \,d\theta ,
\end{align*}
which agrees with the right-hand side of \eqref{;q1plv} since 
\begin{align*}
 \frac{1}{2}\frac{\pr (R_{t}\in dr)}{dr}=\frac{1}{\sqrt{2\pi t^{3}}}r^{2}\exp \left( -\frac{r^{2}}{2t}\right) 
\end{align*}
for $r>0$ (see e.g., \cite[p.\ 252]{ry}) and  
\begin{align*}
 \int _{0}^{\pi }f(\cos \theta )\sin \theta \,d\theta =2\ex [f(V)]
\end{align*}
by a simple change of variables. This ends the proof of the lemma thanks to the arbitrariness of 
$F$ and $f$.
\end{proof}

It is well known (see, e.g., \cite{pr}) that, given three independent copies $N_{1},N_{2},N_{3}$ of 
a nondegenerate, one-dimensional centered Gaussian random variable, the random vector 
\begin{align*}
 \frac{(N_{1},N_{2},N_{3})}{|(N_{1},N_{2},N_{3})|}
\end{align*}
is uniformly distributed on the surface of the two-dimensional unit sphere and 
\begin{align*}
 \frac{N_{1}}{|(N_{1},N_{2},N_{3})|}\eqd V.
\end{align*}
Aside from the problematic fact that 
\begin{align*}
 \int _{0}^{s}\frac{du}{|\bd _{u}|^{2}}=\infty ,\quad 0<s<\infty ,
\end{align*}
a.s., the conclusion~\eqref{;qlv} of the lemma may be inferred from the above fact and 
the skew-product representation of multidimensional Brownian motion, in a stronger 
statement that 
\begin{align}\label{;qlvd}\tag{\ref{;qlv}$'$}
 \left( 
 \{ R_{s}\} _{s\ge 0},VR_{t}
 \right) \eqd 
 \left( 
 \{ |\bd _{s}|\} _{s\ge 0},B^{1}_{t}
 \right) ;
\end{align}
for the skew-product representation, see, e.g., \cite[p.\ 270]{im}. In fact, thanks to the 
Markov property, reasoning of the above proof allows us to obtain \eqref{;qlvd} as well.

Using the above two lemmas, we prove \tsref{;tmaini} and \ref{;tmainii}. 
We start with 

\begin{proof}[Proof of \tref{;tmaini}]
By Pitman's identity~\eqref{;pit}, we see in particular that the process 
\begin{align*}
 \bigl( \cp (B)(s),B_{s}\bigr) ,\quad 0\le s\le t,
\end{align*}
has the same law as 
\begin{align*}
 \Bigl( 
 R_{s},-R_{s}+\min \Bigl\{ 2\min _{s\le u\le t}R_{u}, 2\inf _{u\ge t}R_{u}\Bigr\} 
 \Bigr) ,\quad 0\le s\le t,
\end{align*}
which, by \eqref{;q1lu}, is identical in law with 
\begin{align*}
 \Bigl( 
 R_{s},-R_{s}+\min \Bigl\{ 2\min _{s\le u\le t}R_{u}, R_{t}+(2U-1)R_{t}\Bigr\} 
 \Bigr) ,\quad 0\le s\le t,
\end{align*}
where $U$ is independent of $R$. By noting that $2U-1\eqd V$, the 
conclusion follows thanks to \lref{;lv}.
\end{proof}

The proof of \tref{;tmainii} proceeds similarly.

\begin{proof}[Proof of \tref{;tmainii}]
L\'evy's identity~\eqref{;lev} together with Pitman's identity~\eqref{;pit} entails that 
the process 
\begin{align*}
 \bigl( |B_{s}|+L_{s},|B_{s}|\bigr) ,\quad 0\le s\le t,
\end{align*}
has the same law as 
\begin{align*}
 \Bigl( 
 R_{s},R_{s}-\inf _{u\ge s}R_{u}
 \Bigr) ,\quad 0\le s\le t,
\end{align*}
which, by writing 
\begin{align*}
 \inf _{u\ge s}R_{u}=\min \Bigl\{ 
\min _{s\le u\le t}R_{u},\inf _{u\ge t}R_{u}
\Bigr\} ,\quad 0\le s\le t,
\end{align*}
as in the previous proof, is identical in law with 
\begin{align*}
 \Bigl( 
 R_{s},R_{s}-\min \Bigl\{ 
\min _{s\le u\le t}R_{u},R_{t}-UR_{t}
\Bigr\} 
 \Bigr) ,\quad 0\le s\le t,
\end{align*}
in view of \eqref{;q2lu}. Here, $R$ and $U$ are independent. Since 
\begin{align*}
 \left( 
 \{ R_{s}\} _{0\le s\le t},UR_{t}
 \right) \eqd 
 \left( 
 \{ |\bd _{s}|\} _{0\le s\le t},|B^{1}_{t}|
 \right) 
\end{align*}
by \eqref{;absv} and \lref{;lv}, we have the theorem.
\end{proof}

\begin{remark}
If we apply \eqref{;q1lu} instead of \eqref{;q2lu}, then we obtain another  
equivalence in law: the process~\eqref{;q1tmainii} has the same law as 
\begin{align*}
 \Bigl( 
 |\bd _{s}|,|\bd _{s}|-\min \Bigl\{ 
 \min _{s\le u\le t}|\bd _{u}|,|B^{1}_{t}|
 \Bigr\} 
 \Bigr) ,\quad 0\le s\le t,
\end{align*}
from which \cref{;cmainii} does not immediately follow. 
\end{remark}

\section{Proofs of corollaries}\label{;spcs}

In this section, we provide proofs of corollaries to the two 
\tsref{;tmaini} and \ref{;tmainii} 
proven in the previous section. 

\begin{proof}[Proof of \cref{;cmaini}]
The assertion follows readily from \eqref{;jil} because of the independence of 
$B^{1}$ and $(B^{2},B^{3})$, and the definition of the process $M^{x}$ given 
just below \tref{;tmaini}.
\end{proof}

\begin{proof}[Proof of \cref{;ccmaini}]
The process $M^{x}$ is identical in law with $\cp (\bb ^{x})$ 
by \cref{;cmaini}. Taking $y=0$ therein, we know from \cite[Theorem~1.4]{har25+} 
that the law of $\cp (\bb ^{x})$ agrees with that of $\cp (\bb ^{0})$ 
conditioned on the event that $\cp (\bb ^{0})(t)\ge |x|$ (see also \rref{;rce} below), 
which, by \cref{;cmaini} again, proves the assertion.
\end{proof}

\begin{remark}\label{;rce}
From the proof of \tref{;tmaini} in the previous section, we see that 
\begin{align*}
 \bigl( 
 \{ \cp (B)(s)\} _{0\le s\le t},B_{t}
 \bigr) 
 \eqd 
 \bigl( 
 \{ R_{s}\} _{0\le s\le t},VR_{t}
 \bigr) ,
\end{align*}
where $V$ is independent of $R$. This identity in law implies that 
\begin{align*}
 \ex \!\left[ 
 F\bigl( \cp (\bb ^{x})\bigr) 
 \right] \frac{1}{\sqrt{2\pi t}}\exp \left( 
 -\frac{x^{2}}{2t}
 \right) =\ex \!\left[ 
 F(R)\frac{1}{2R_{t}};R_{t}\ge |x|
 \right] ,
\end{align*}
where $F$ is any nonnegative measurable function on $\ctd{t}{}$ and 
$x\in \R $ is arbitrary. Comparison with the case $x=0$ leads to 
the conditional equivalence between $\cp (\bb ^{x})$ and $\cp (\bb ^{0})$ 
mentioned above. In addition, from Imhof's relation~\eqref{;imrel} it follows that 
\begin{align*}
 \ex \!\left[ 
 F\bigl( \cp (\bb ^{x})\bigr) 
 \right] \exp \left( 
 -\frac{x^{2}}{2t}
 \right) =\ex \!\left[ 
 F(M^{0});M^{0}_{t}\ge |x|
 \right] .
\end{align*}
\end{remark}

We proceed to the proof of \cref{;cmainii}.

\begin{proof}[Proof of \cref{;cmainii}]
As mentioned earlier, the process $\{ B_{s}\} _{0\le s\le t}$ given $|B_{t}|=|x|$ has 
the same law as $\{ \ve \bb ^{x}_{s}\} _{0\le s\le t}$, with $\ve $ a Rademacher 
random variable independent of $\bb ^{x}$. This conditional equivalence together with 
\eqref{;chard} shows that, conditionally on $|B_{t}|=|x|$, the law of the process 
\eqref{;q1tmainii} agrees with that of \eqref{;q1cmainii}. Moreover, the above 
equivalence in law applied to $\{ B^{1}_{s}\} _{0\le s\le t}$ entails that, thanks to 
the independence of $B^{1}$ and $(B^{2},B^{3})$, the process $\{ |\bd _{s}|\} _{0\le s\le t}$ 
given $|B^{1}_{t}|=|x|$ is identical in law with 
\begin{align*}
 \sqrt{(\ve \bb ^{x}_{s})^{2}+(B^{2}_{s})^{2}+(B^{3}_{s})^{2}},\quad 0\le s\le t,
\end{align*}
and hence with $M^{x}$ as $\ve ^{2}=1$ a.s. Here, $\ve ,\bb ^{x}$ and 
$(B^{2},B^{3})$ are independent. Consequently, the law of the process~\eqref{;q2tmainii} 
given $|B^{1}_{t}|=|x|$ coincides with that of \eqref{;q2cmainii}, which completes the 
proof. 
\end{proof}

We conclude this paper with the proof of \cref{;ccmainii}. In view of \cref{;cmainii} and 
\eqref{;eqds}, it suffices to show that 
\begin{equation}\label{;suff}
\begin{split}
 &2\max _{0\le u\le s}\phi _{u}-\min \Bigl\{ 
 \min _{s\le u\le t}\cp (\phi )(u),\cp (\phi )(t)-|x|
 \Bigr\} \\
 &=\min \Bigl\{ 
 \max _{0\le u\le s}\phi _{u},\Bigl( \max _{s\le u\le t}\phi _{u}\Bigr) _{+}
 \Bigr\} ,
\end{split}
\end{equation}
for all $0\le s\le t$ and for all $\phi \in \ctd{t}{}$ such that 
\begin{align}\label{;bc}
 \phi _{0}=0 && \text{and} && \phi _{t}=-|x|.
\end{align}
The following lemma is a restatement of \cite[Lemma~2.2]{har25+}.

\begin{lemma}[{{\cite{har25+}}}]\label{;lpit}
For every $\phi \in C([0,t];\R )$, it holds that 
\begin{align*}
 \min _{s\le u\le t}\cp (\phi )(u)
 =2\max _{0\le u\le s}\phi _{u}-\min \Bigl\{ 
 \max _{0\le u\le s}\phi _{u},\max _{s\le u\le t}\phi _{u}
 \Bigr\} 
\end{align*}
for all $0\le s\le t$.
\end{lemma}

\begin{remark}\label{;rscl}
By the above lemma, we see that, for every $0\le s\le t$, the left-hand side of 
\eqref{;scl} is equal to 
\begin{align*}
 \phi _{s}+2\min \Bigl\{ 
 \Bigl( \max _{0\le u\le s}\phi _{u}-\max _{s\le u\le t}\phi _{u}\Bigr) _{+},
 \max _{0\le u\le t}\phi _{u}-\max _{0\le u\le s}\phi _{u}
 \Bigr\} ,
\end{align*}
which agrees with $\phi _{s}$ since 
\begin{align*}
 \max _{0\le u\le t}\phi _{u}-\max _{0\le u\le s}\phi _{u}
 =\Bigl( \max _{s\le u\le t}\phi _{u}-\max _{0\le u\le s}\phi _{u}\Bigr) _{+}.
\end{align*}
\end{remark}

By using \lref{;lpit}, the proof of \cref{;ccmainii} proceeds as follows: 

\begin{proof}[Proof of \cref{;ccmainii}]
Fix $s\in [0,t]$ and pick arbitrarily $\phi \in \ctd{t}{}$ satisfying 
\eqref{;bc}. Then, by \lref{;lpit} and by the latter condition in \eqref{;bc}, 
the left-hand side of \eqref{;suff} is written as 
\begin{align}\label{;llhs}
 &2\max _{0\le u\le s}\phi _{u}-\min \Bigl\{ 
 2\max _{0\le u\le s}\phi _{u}-\min \Bigl\{ 
 \max _{0\le u\le s}\phi _{u},\max _{s\le u\le t}\phi _{u}
 \Bigr\} ,2\max _{0\le u\le t}\phi _{u}
 \Bigr\} \notag \\
 &=-\min \bigl\{ 
 -\min \{ a,b\}, 2(\max \{ a,b\} -a)  
 \bigr\} \notag \\
 &=\max \bigl\{ 
 \min \{ a,b\} ,2(a-\max \{ a,b\} )
 \bigr\}  ,
\end{align}
where we put 
\begin{align*}
 a=\max _{0\le u\le s}\phi _{u}, && b=\max _{s\le u\le t}\phi _{u}.
\end{align*}
Notice that $a\ge 0$ by the former condition in \eqref{;bc}. Therefore, when 
$a<b$, the last expression~\eqref{;llhs} of the left-hand side of 
\eqref{;suff} is equal to 
\begin{align*}
 \max \{ a, 2(a-b)\} =a 
\end{align*}
for $2(a-b)<0$, which agrees with the right-hand side of \eqref{;suff}; indeed,  
\begin{align*}
 \min \{ a,b_{+}\} &=\min \{ a,b\} \\
 &=a,
\end{align*}
where, for the first line, we have used the positivity of $b$: $b>a\ge 0$.
On the other hand, when $a\ge b$, the expression~\eqref{;llhs} is equal to 
\begin{align*}
 \max \{ b,0\} =b_{+},
\end{align*}
which agrees with the right-hand side of \eqref{;suff} because 
\begin{align*}
 \min \{ a,b_{+}\} =b_{+}
\end{align*}
due to the fact that $a\ge b$ and $a\ge 0$. This completes the proof 
of \eqref{;suff} and \cref{;ccmainii} is proven.
\end{proof}

\medskip 
\noindent 
{\bf Funding} This work was supported in part by JSPS KAKENHI Grant Number  22K03330.



\begin{thebibliography}{99}

\bibitem{ber0} J.~Bertoin, D\'ecomposition du mouvement brownien avec d\'erive en un minimum local par juxtaposition de ses excursions positives et n\'egatives, in: S\'eminaire de Probabilit\'es, XXV, pp.~330--344, Lecture Notes in Math.\ {\bf 1485}, Springer, Berlin, 1991.

\bibitem{ber} J.~Bertoin, An extension of Pitman's theorem for spectrally positive L\'evy processes, Ann.\ Probab.\ {\bf 20} (1992), 1464--1483.

\bibitem{bp} J.~Bertoin, J.~Pitman, Path transformations connecting Brownian bridge, excursion and meander, Bull.\ Sci.\ Math.\ {\bf 118} (1994), 147--166.

\bibitem{by} P.~Biane, M.~Yor, Quelques pr\'ecisions sur le m\'eandre brownien, Bull.\ Sci.\ Math.\ 2e S\'erie {\bf 112} (1988), 101--109.

\bibitem{har25+} Y.~Hariya, Invariance of three-dimensional Bessel bridges in terms of time reversal, Electron.\ J.\ Probab.\ {\bf 31} (2026),  Paper No.~92, 24 pp.

\bibitem{imh} J.-P.~Imhof, Density factorizations for Brownian motion, meander and the three-dimensional Bessel process, and applications, J.\ Appl.\ Probab.\ {\bf 21} (1984), 500--510.

\bibitem{im}  K.~It\^o, H.P.~McKean, Jr., Diffusion Processes and their Sample Paths, second printing, corrected, Springer, Berlin, 1974.

\bibitem{lg} J.-F.~Le Gall, Brownian Motion, Martingales, and Stochastic Calculus, Springer, Cham, 2016.

\bibitem{my} R.~Mansuy, M.~Yor, Aspects of Brownian Motion, Springer, Berlin, 2008.

\bibitem{mySII} H.~Matsumoto, M.~Yor, Exponential functionals of Brownian motion, II: Some related diffusion processes, Probab.\ Surv.\ {\bf 2} (2005), 348--384. 

\bibitem{jwp} J.W.~Pitman, One-dimensional Brownian motion and the three-dimensional Bessel process, Adv.\ in Appl.\ Probab.\ {\bf 7} (1975), 511--526.

\bibitem{pr} J.~Pitman, N.~Ross, Archimedes, Gauss, and Stein, Notices Amer.\ Math.\ Soc.\ {\bf 59} (2012), 1416--1421.

\bibitem{ry} D.~Revuz, M.~Yor, Continuous Martingales and Brownian Motion, 3rd ed., Springer, Berlin, 1999.

\bibitem{rv} B.~Rider, B.~Valk\'o, Matrix Dufresne identities, Int.\ Math.\ Res.\ Not.\ {\bf 2016} (2016), 174--218.

\bibitem{rvy} B.~Roynette, P.~Vallois, M.~Yor, Some extensions of Pitman and Ray--Knight theorems for penalized Brownian motions and their local times, IV, Studia Sci.\ Math.\ Hungar.\ {\bf 44} (2007), 469--516.

\bibitem{st} Y.~Saisho, H.~Tanemura, Pitman type theorem for one-dimensional diffusion processes, Tokyo J.\ Math.\ {\bf 13} (1990), 429--440.

\end{thebibliography}
\end{document}